\documentclass[10pt]{article}

\usepackage[utf8]{inputenc}
\usepackage{amsthm}
\usepackage{amssymb}
\usepackage{mathtools}
\usepackage{graphicx}	
\usepackage{pgf}
\usepackage{tikz,pgfplots}
\usepackage{extarrows}
\pgfplotsset{compat=1.15}
\usepackage{booktabs}
\usepackage[pdftex,
              bookmarks,
              bookmarksnumbered=true,
              bookmarksopen=false,
              colorlinks,
              citecolor=black,
              filecolor=black,
              linkcolor=black,
              urlcolor=black,
              pdfauthor={Alejandro Garces},
              pdftitle={A review of matrix algebra for power and energy applications},
              pdftoolbar=true,
              pdfmenubar=true,
              pdffitwindow = true
              ]{hyperref}
\usepackage{xcolor}
\selectcolormodel{gray}

\newcommand{\julia}{
\begin{tikzpicture}[x=1mm,y=1mm]
\fill[gray] (0,0) circle (0.9);
\fill[gray] (2,0) circle (0.9);
\fill[gray] (1,1.6) circle (0.9);
\end{tikzpicture}
}
\newcounter{example}
\newenvironment{example}{\refstepcounter{example}\par\medskip
   \noindent\textbf{\julia Example~\theexample.  \\ } \rmfamily}{\hfill \\ 
   \textcolor{gray}{\rule{15cm}{1pt}}\medskip}

\usepackage{listings}

\lstdefinelanguage{Julia}%
  {morekeywords={abstract,break,case,catch,const,continue,do,else,elseif,%
      end,export,false,for,function,immutable,import,importall,if,in,%
      macro,module,otherwise,quote,return,switch,true,try,type,typealias,%
      using,while},%
   sensitive=true,%
   alsoother={$.$},%
   morecomment=[l]\#,%
   morecomment=[n]{\#=}{=\#},%
   morestring=[s]{"}{"},%
   morestring=[m]{'}{'},%
}[keywords,comments,strings]%
\usepackage[top=3cm,bottom=2.5cm,right=2.5cm,left=2.5cm]{geometry}

\usepackage{stix2}
\usepackage[font={bf,small}]{caption}

\usepackage{fancyhdr}
\begin{document}
\title{A Review of Matrix Algebra for Power and Energy Applications}
\author{Alejandro Garcés-Ruiz}
\date{Universidad Tecnológica de Pereira \\ Carrera 27 nro. 10-02 Los Álamos-Pereria-Risaralda-Colombia \\ www.utp.edu.co}

\maketitle

\abstract{This report presents a brief review of matrix algebra and its implementation in Julia for power and energy applications. First, we present basic examples of data visualization, followed by conventional operations with matrices and vectors.  Then, we study quadratic forms and norms, two main concepts required in the convergence study of the power flow in power and energy applications.  After that, we give good practices to create a neat code in Julia. There is an extensive set of examples available on the Internet related to these basic aspects, so we avoid repeating what is well documented.  Hence, we show only basic examples to create our first scripts.}

\tableofcontents

\section{An introduction to Julia}

Julia is a high-level, high-performance programming language designed for technical computing. It offers several advantages, which have contributed to its growing popularity in power systems applications. The syntax is easy to learn and read, since it resembles other high-level languages like Matlab and Python, two of the most common languages for programming computational methods in power systems and power distribution networks. However, Julia outperforms these languages, approaching the performance of low-level languages like C and Fortran.

Like Python, it is an open-source project with an active community of developers and researchers who create packages, and provide support through forums, web pages, and mailing lists. There are different environments for developing Julia code such as Visual Studio Code, Juno and Jupyter Notebooks. We do not catch ourselves with any of these environments, since Julia is a constantly evolving ecosystem.  We recommend checking the latest information on Julia's official website and the community resources for the most up-to-date details.  Let's introduce our first examples:

\begin{example}
Julia has a built-in package manager that ensures seamless integration with external libraries and tools.  A new package is added using the package manager as follows:
\begin{lstlisting}[frame=trBL]
using Pkg
Pkg.add("Plots")
\end{lstlisting}
This code installs the package \texttt{Plots} which is an extensive library for data visualization. The next example shows its usage.  
\end{example}

\begin{example} 
\texttt{Plots} is a library that resembles Matlab functionalities. Once the package is installed, we can plot a simple function as given below:
\begin{lstlisting}[frame=trBL]
using Plots
x = LinRange(0,2*pi,100)
y = sin.(x)
plot(x,y)
\end{lstlisting}
The first line calls the package \texttt{Plots} that allows data visualization, while the second line generates an array with 100 points between $0$ and $2\pi$. Julia allows special characters into the code, so we can replace the word \texttt{pi} by the Greek letter \texttt{$\pi$}, obtaining the same result. The third line applies the function \texttt{sin} to each point $x$. Finally, the fourth line plots $x$ vs $y$. The reader is invited to implement this code to see the result. Notice the function \texttt{sin} is evaluated in each point of the array \texttt{x}. Hence, it requires a dot to indicate it is a pointwise operation.
\end{example}

\begin{example}
Julia resembles both Python and Matlab; thus, it is easy to translate a code from any of these languages to Julia. Let us see two ways to define a for loop and print the iteration:
\begin{lstlisting}[frame=trBL]
for k = 1:10
    println("iteration:",k)
end

for k in 1:10
     println("iteration:",k)
end
\end{lstlisting}
Both codes are functionally equivalent, and the choice between them depends on the specific application and personal preferences. Readers familiar with Matlab may prefer the first code, whereas those accustomed to Python may favor the second.
\end{example}
\begin{example}
We can also define a loop inside an array in order to define its elements, namely:
\begin{lstlisting}[frame=trBL]
using Plots
v = [cos(t)+sin(t)*1im for t in 0:0.1:2*pi]
plot(v)
\end{lstlisting}
Notice the imaginary number is given by \texttt{1im}. In this case, we have created a complex vector $v = \cos(t)+i\sin(t)$. Then, we plot the real and imaginary parts, obtaining a unitary circle. \index{im}\index{complex numbers}
\end{example}

\section{Vectors, matrices and hyper-arrays}\label{sec:propiedades_traza_determinante}

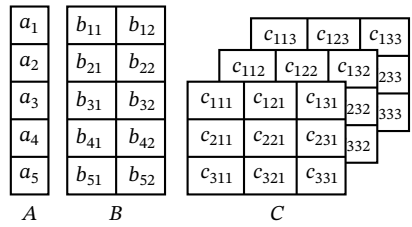
\begin{figure}[hbt]
\centering
\footnotesize
\begin{tikzpicture}[x=1mm,y=1mm]
	\draw[thick] (-2.5,2.5) rectangle +(5,25);
	\draw[thick] (-2.5,7.5) -- +(5,0);
	\draw[thick] (-2.5,12.5) -- +(5,0);
	\draw[thick] (-2.5,17.5) -- +(5,0);
	\draw[thick] (-2.5,22.5) -- +(5,0);
	\node at (0,25) {$a_1$};
	\node at (0,20) {$a_2$};
	\node at (0,15) {$a_3$};
	\node at (0,10) {$a_4$};
	\node at (0,5) {$a_5$};
	\node at (0,0) {$A$};
	
	\draw[thick] (5,2.5) rectangle +(13,25);
	\draw[thick] (5,7.5) -- +(13,0);
	\draw[thick] (5,12.5) -- +(13,0);
	\draw[thick] (5,17.5) -- +(13,0);
	\draw[thick] (5,22.5) -- +(13,0);
	\draw[thick] (11.5,2.5) -- +(0,25);
	\node at (8,25) {$b_{11}$};
	\node at (8,20) {$b_{21}$};
	\node at (8,15) {$b_{31}$};
	\node at (8,10) {$b_{41}$};
	\node at (8,5) {$b_{51}$};
	\node at (15,25) {$b_{12}$};
	\node at (15,20) {$b_{22}$};
	\node at (15,15) {$b_{32}$};
	\node at (15,10) {$b_{42}$};
	\node at (15,5) {$b_{52}$};
	\node at (11.5,0) {$B$};
	
	\begin{scope}[xshift=24,yshift=24]
		\draw[thick,fill=white] (21,2.5) rectangle +(21,15);
		\draw[thick] (21,7.5) -- +(21,0);
		\draw[thick] (21,12.5) -- +(21,0);
		\draw[thick] (28.5,2.5) -- +(0,15);
		\draw[thick] (35.5,2.5) -- +(0,15);
		\node at (25,15) {$c_{113}$};
		\node at (32,15) {$c_{123}$};
		\node at (39,15) {$c_{133}$};
		\node at (39,10) {$c_{233}$};
		\node at (39,5) {$c_{333}$};
	\end{scope}

	\begin{scope}[xshift=12,yshift=12]
		\draw[thick,fill=white] (21,2.5) rectangle +(21,15);
		\draw[thick] (21,7.5) -- +(21,0);
		\draw[thick] (21,12.5) -- +(21,0);
		\draw[thick] (28.5,2.5) -- +(0,15);
		\draw[thick] (35.5,2.5) -- +(0,15);
		\node at (25,15) {$c_{112}$};
		\node at (32,15) {$c_{122}$};
		\node at (39,15) {$c_{132}$};
		\node at (39,10) {$c_{232}$};
		\node at (39,5) {$c_{332}$};
	\end{scope}
	
	\draw[thick,fill=white] (21,2.5) rectangle +(21,15);
	\draw[thick] (21,7.5) -- +(21,0);
	\draw[thick] (21,12.5) -- +(21,0);
	\draw[thick] (28.5,2.5) -- +(0,15);
	\draw[thick] (35.5,2.5) -- +(0,15);
	\node at (25,15) {$c_{111}$};
	\node at (25,10) {$c_{211}$};
	\node at (25,5) {$c_{311}$};
	\node at (32,15) {$c_{121}$};
	\node at (32,10) {$c_{221}$};
	\node at (32,5) {$c_{321}$};
	\node at (39,15) {$c_{131}$};
	\node at (39,10) {$c_{231}$};
	\node at (39,5) {$c_{331}$};
	\node at (33,0) {$C$};
\end{tikzpicture}
    \caption{Example of a vector, a matrix, and a three-dimensional hyper-array}
    \label{fig:ch2_tensorgrafico}
\end{figure}

Some physical magnitudes, such as velocity or acceleration, are characterized by having not only magnitude but also direction. These magnitudes are called vectors and are usually defined in the three-dimensional space; however, we can generalize the concept to arrays in $\mathbb{R}^n$, $\mathbb{C}^n$. In the same way, we can build two-dimensional arrays, that is, matrices and even three-dimensional arrays (tensors or hyper-arrays) that resemble Rubik's cubes, as shown in Figure \ref{fig:ch2_tensorgrafico}.

These arrays are easily represented in Julia, as presented in the following example:

\begin{example}
This example creates random arrays of different sizes. Entries of these arrays are obtained using square brackets, for instance: \texttt{A[1]}, \texttt{B[4,1]},\texttt{C[3,2,1]}; other outputs are given below:
\begin{lstlisting}[frame=trBL]
A = rand(5)
B = rand(5,2)
C = rand(3,3,3)
println("Vector A =", A)
println("Matrix B =", B)
println("Hyper-array C =", C)
println("First column of B =", B[:,1])
println("Second row of B   =", B[2,:])
println("First matrix of C =", C[:,:,1])
\end{lstlisting}
\end{example}

Most of the problems in power systems are represented in high dimensional spaces.  Therefore, we require storing and operate with vectors, matrices and hyper-arrays. Basic operations such as adding and multiplying by a scalar are intuitively defined for both, real and complex domains.

\begin{example}
This example shows basic operations with complex arrays:
\begin{lstlisting}[frame=trBL]
A = [3+2im 8+3im]
B = [1-4im 7+5im]
C = A+B
println(C)
println(C')
println(transpose(C))
\end{lstlisting}
Notice the difference between \texttt{C'} and \texttt{transpose(C)}.  The former returns the transpose conjugate, whereas the latter returns the transpose without conjugate.
\end{example}

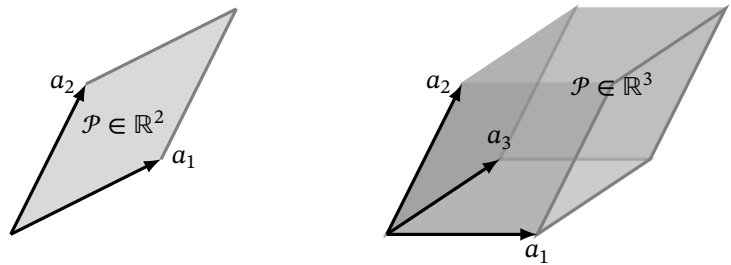
\begin{figure}[htb]
    \centering
    \begin{tikzpicture}[x=1mm,y=1mm]
        \fill[gray!30] (0,0) -- (20,10) -- (30,30) -- (10,20) -- cycle;
        \draw[very thick,gray] (10,20) -- +(20,10);
        \draw[very thick,gray] (20,10) -- +(10,20);
        \draw[-latex,very thick] (0,0) -- (20,10) node[right] {$a_1$};
        \draw[-latex,very thick] (0,0) -- (10,20) node[left] {$a_2$};
        \node at (15,15) {$\mathcal{P}\in\mathbb{R}^2$};
    \begin{scope}[shift={(50,0)}]
        \draw[very thick,gray,fill,opacity=0.4] 
        (0,0) -- (15,10) -- (35,10) -- (20,0) -- cycle;
        \draw[very thick,gray,fill,opacity=0.6] 
        (0,0) -- (10,20) -- (25,30) -- (15,10) -- cycle;
        \draw[very thick,gray,fill,opacity=0.5] 
        (15,10) -- (25,30) -- (45,30) -- (35,10) -- cycle;
        \draw[very thick,gray] 
        (20,0) -- (35,10) -- (45,30) -- (30,20) -- cycle;
        \draw[very thick,gray,fill,opacity=0.3] 
        (0,0) -- (20,0) -- (30,20) -- (10,20) -- cycle;
        
        \draw[-latex,very thick] (0,0) -- (20,0) node[below] {$a_1$};
        \draw[-latex,very thick] (0,0) -- (10,20) node[left] {$a_2$};
        \draw[-latex,very thick] (0,0) -- (15,10) node[above] {$a_3$};
        \node at (30,20) {$\mathcal{P}\in\mathbb{R}^3$};
    \end{scope}            
    \end{tikzpicture}
    \caption{Example of a parallelogram ($\mathbb{R}^2$) or a parallepiped ($\mathbb{R}^3$) formed by a set of linearly independent vectors. }
    \label{fig:paralelepipedo}
\end{figure}

We can define other operators such the trace, determinant, and transpose for matrix variables. Let us revisit briefly these concepts:

The trace of a square matrix $A$ is denoted by $\operatorname{tr}(A)$. It is a function $\operatorname{tr}:\mathbb{R}^{n\times n}\rightarrow \mathbb{R}$ that takes the sum of the elements in the diagonal, namely:
\begin{equation}
\operatorname{tr}(A) = \sum_{k=1}^{n} a_{kk}
\end{equation}
The following properties hold for any pair of square matrices $A,B$ and, a scalar $\alpha$:
\begin{align}
    \operatorname{tr}(A+B) &=\operatorname{tr}(A) + \operatorname{tr}(B) \\
    \operatorname{tr}(\alpha A) &= \alpha \operatorname{tr}(A) \\
    \operatorname{tr}\left(A^\top\right) &= \operatorname{tr}(A) \\    
    \left(A^\top\right)^\top &= A \\
    (A+B)^\top &= A^\top + B^\top \\
    (AB)^\top &= B^\top A^\top
\end{align}
The determinant is a another matrix function $\operatorname{det}:\mathbb{R}^{n\times n}\rightarrow \mathbb{R}$, that takes the $n$ column (or row) vectors that form a matrix, and returns a real number $\operatorname{det}(A)=\operatorname{det}([a_1,a_2,\dots,a_n])$. There are several interpretations of the determinant. Perhaps the most useful, and visually meaningful of these interpretations, is the volume of the parallelogram generated by the vectors that form the matrix (see Figure \ref{fig:paralelepipedo}). 

The area of the parallelogram $\mathcal{P}$ formed by the vectors $v_1$ and $v_2$ can be calculated by organizing these vectors into a matrix $A$ as given below:
\begin{equation}
    \operatorname{Area}(\mathcal{P}) = |\operatorname{det}(A)|
\end{equation}
This concept can be extended to the volume of $\mathcal{P}$ for the three-dimensional case and the hyper-volume for higher dimensions.

Interpreting the determinant as a volume helps to remember and better understand some of its properties.  For example, $\operatorname{det}(A)=0$ if $a_1$ and $a_2$ are linearly dependent because in that case, both vectors are parallel (or antiparallel) and do not generate any area.  The determinant of the identity matrix is $\operatorname{det}(I)=1$ because it represents a square generated by two vectors of size one. Other properties are given below:
\begin{align}
\operatorname{det}\left(A^\top\right) &= \operatorname{det}(A) \\
\operatorname{det}\left(A^{-1}\right) &= 1/\operatorname{det}(A) \\
\operatorname{det}(AB) &= \operatorname{det}(A)\operatorname{det}(B) \\
\operatorname{det}\left(P^{-1}AP\right) &= \operatorname{det}(A) \\
\operatorname{det}(A) &\neq 0, \; \text{if } A \text{ is non-singular}
\end{align}

\begin{example}
The package \texttt{LinearAlgebra} allows operating with vectors and matrices.  The following code defines a $2\times 2$ matrix of zeros and calculates its trace, determinant, and rank:
\begin{lstlisting}[frame=trBL]
using LinearAlgebra     
M = zeros(2,2)          
d = det(M)              
t = tr(M)               
r = rank(M)             
printstyled("determinant of M \t = ",
d, "\n", color = :red) 
printstyled("trace of M \t = ",t, "\n",
color = :green) 
printstyled("rank of M  \t = ",r, "\n", 
color = :blue) 
\end{lstlisting}
Additional features of Julia were showcased in this script. Comments are defined with the number symbol (hashtag); character \texttt{$\setminus$t} inserts a tap, while \texttt{$\setminus$n} creates a new line; \texttt{printstyled} allows printing the results with different colors and format.
\end{example}

Besides the conventional matrix multiplication, we can also define the Hadamard and Kronecker product.  Given two matrices $A$ and $B$, the Hadamard product, denoted by $A\circ B$, is the obtained by an element-wise multiplication; for this, the size of both matrices must be the same. The Kronecker product returns an augmented matrix, denoted by $A\otimes B$, where each element of $A$ is multiplied by the entire matrix $B$.

\begin{example}
Consider two matrices $A$ and $B$ defined as follows:
\begin{align}
A &= \begin{bmatrix} 1 & 2 \\ 8 & 5 \end{bmatrix} \\
B &= \begin{bmatrix} 3 & 4 \\ 7 & 6 \end{bmatrix}
\end{align}
The following code returns the conventional product, the Hadamard product and, the Kronecker product:
\begin{lstlisting}[frame=trBL]
using LinearAlgebra
A = [1 2; 8 5]
B = [3 4;7 6]
println("Conventional product\t =", A*B)
println("Hadamard product\t =", A.*B)
println("Kronecker product\t =", kron(A,B))
\end{lstlisting}
For the sake of better understanding the concept, we present each of these operations step by step:
\begin{align}
AB &= \begin{bmatrix}
(1)(3)+(2)(7) & (1)(4)+(2)(6) \\ 
(8)(3)+(5)(7) & (8)(4)+(5)(6) \end{bmatrix} =
\begin{bmatrix} 17 & 16 \\ 59 & 62 \end{bmatrix}
\\
A\circ B &= \begin{bmatrix}
(1)(3) & (2)(4) \\ 
(8)(7) & (5)(6) \end{bmatrix}
= \begin{bmatrix} 3 & 8 \\ 56 & 30 \end{bmatrix}
\\
A\otimes B &= \begin{bmatrix} 
(1)\begin{bmatrix} 3 & 4 \\ 7 & 6 \end{bmatrix}
&
(2)\begin{bmatrix} 3 & 4 \\ 7 & 6 \end{bmatrix}
& \\
\\
(8)\begin{bmatrix} 3 & 4 \\ 7 & 6 \end{bmatrix}
&
(5)\begin{bmatrix} 3 & 4 \\ 7 & 6 \end{bmatrix}
\end{bmatrix}
=\begin{bmatrix} 
3 & 4 & 6 & 8 \\ 7 & 6 & 14 & 12\\ 24 & 32 & 15 & 20\\  56 & 48 & 35 & 30
\end{bmatrix}
\end{align}
Very complicated operations in power systems can be represented and analyzed by a Hadamard or a Kronecker product.
\end{example}

It is important to identify the properties of each of these products. For instance, conventional matrix multiplication requires that the number of columns of the first matrix is the same as the number of rows of the second matrix. This multiplication is usually non-commutative ($AB\neq BA$).  The Hadamard product requires two matrices of the same size.  
This product is commutative, i.e., $A\circ B = B\circ A$. The Kronecker product allows matrices of any size.  Its result is non-commutative.

\section{Linear transformations, eigenvalues, and eigenvectors}

A linear transformation is a mathematical operation that takes input points and transforms them, creating another space where, perhaps, it is easy to make calculations.  This transformation involves stretching, squeezing, rotating, reflecting, or shearing the points. However, it must follow certain rules such as preserving the origin, and maintaining addition and scalar multiplication.   Linear transformations appear in several power systems' computation problems.  These include the Park and Clarke transformation as well as the symmetrical components.

Matrix transformation are usually represented by a transformation matrix $A$.  Thus, each vector $x$ from the original space can be transformed into a new vector $z = Ax$. Geometric properties of the transformed space are given by $A$. There exist some vectors $v$ which maintain their direction in the original and the transformed space.  That is to say, $Av = \lambda v$.  These are the eigenvectors of $A$ and the proportional constant $\lambda$ are the eigenvalues.  Every square matrix has a set of unique eigenvalues. However, eigenvectors are not unique, since every multiple is also an eigenvector. Besides, the following two properties hold:
\begin{align}
    \det(A) &= \prod_{i=1}^n \lambda_i \\
    \operatorname{tr}(A) &= \sum_{i=1}^n \lambda_i
\end{align}
Both, eigenvectors and eigenvalues play a key role in stability problems for power systems applications.  Let us consider an example where both concepts are presented in the context of a linear transformation.

\begin{example}
Consider a set of points representing a unit circle in the plane:
\begin{equation}
\Omega = \left\{x \in\mathbb{R}^2: x_1^2+x_2^2 = 1  \right\}
\end{equation}
This set can be created parametrically as given below:
\begin{lstlisting}[frame=trBL]
using LinearAlgebra
using Plots
Omega = zeros(2,100)
for k = 1:100
    Omega[1,k] = cos(2*pi*k/100)
    Omega[2,k] = sin(2*pi*k/100)
end
p1 = plot(Omega[1,:],Omega[2,:],aspect_ratio=:equal)
\end{lstlisting}
Consider now a linear transformation represented by a $2\times 2$ matrix $A$, and a new space $\Gamma$ as follows:
\begin{equation}
\Gamma = \left\{z \in\mathbb{R}^2: z=Ax, \forall x\in\Omega  \right\}
\end{equation}
This set can also be defined and plotted in Julia:
\begin{lstlisting}[frame=trBL]
A = [1 0.3;0.5 1]
Gamma = zeros(2,100)
for k = 1:100
    Gamma[:,k] = A*Omega[:,k]
end
p2 = plot(Gamma[1,:],Gamma[2,:],aspect_ratio=:equal)
\end{lstlisting}
Finally, we can calculate the eigenvalues and eigenvectors of the matrix transformation.  The eigenvectors are the directions that maintain unchanged after the transformation. 
\begin{lstlisting}[frame=trBL]
(L,V) = eigen(A)
p2=plot!([0,V[1,1]],[0,V[1,2]], arrow=1)
p2=plot!([0,V[2,1]],[0,V[2,2]], arrow=1)
plot(p1,p2,layout=(1,2))
\end{lstlisting}

Figure \ref{fig:transformacion_lineal} shows the original and the transformed space. The former contains a circle that is transformed into an ellipse.  The eigenvalues are the major axis of this ellipse. 

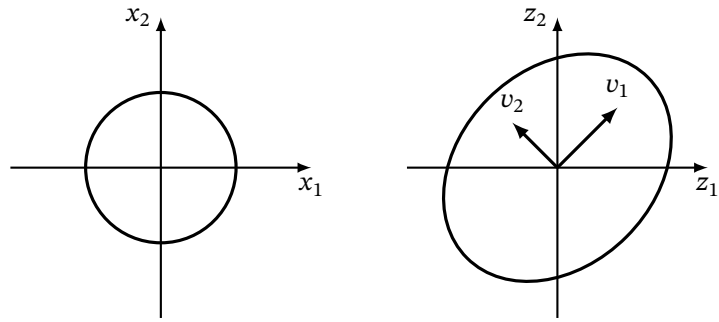
\begin{figure}[hbt]
\centering
\begin{tikzpicture}[x=1mm,y=1mm,thick]
\draw[-latex] (-20,0) -- (20,0) node[below] {$x_1$};
\draw[-latex] (0,-20) -- (0,20) node[left] {$x_2$};
\draw[very thick] (0,0) circle (10);

\begin{scope}[xshift = 150]
\draw[-latex] (-20,0) -- (20,0) node[below] {$z_1$};
\draw[-latex] (0,-20) -- (0,20) node[left] {$z_2$};
\draw[very thick,rotate=-45] (0,0) ellipse (13 and 17);
\draw[-latex, very thick] (0,0) -- (8,8) node[above] {$v_1$};
\draw[-latex, very thick] (0,0) -- (-6,6) node[above] {$v_2$};

\end{scope}
    
\end{tikzpicture}
\caption{Example of a linear transformation.  The eigenvectors $v$ represent the direction which remains after the transformation.}
\label{fig:transformacion_lineal}
\end{figure}
\end{example}

\section{Quadratic forms, semidefinite matrices and norms}\label{sec:cuadraticformss}

A quadratic form is a function from $\mathbb{R}^n$ to $\mathbb{R}$ defined by the following matrix equation:
\begin{equation}
q(x) = x^\top M x \label{eq:quadratic_form}
\end{equation}
where $M$ is a square matrix. This type of function is a polynomial where the maximum power is two.  For example, the following quadratic equation:
\begin{equation}
q(x_1,x_2) = x_1^2 + 3x_1x_2 - x_2^2
\end{equation}
can be represented by a function as \eqref{eq:quadratic_form} with,
\begin{equation}
M = \begin{bmatrix} 1 & 3 \\ 0 & -1\end{bmatrix} \label{eq:MMMM}
\end{equation}

A quadratic form $q$ can be represented by different matrices $M$. For example, we can replace $M$ by a symmetric matrix, $(M^\top + M)/2$ obtaining the same polynomial, namely:
\begin{equation}
q(x_1,x_2) = \begin{bmatrix}x_1 \\ x_2 \end{bmatrix}^\top
\begin{bmatrix} 1 &  {3}/{2} \\ {3}/{2} & -1 \end{bmatrix}
\begin{bmatrix}x_1 \\ x_2 \end{bmatrix}
\end{equation}

Quadratic functions appear in many applications in power systems, among which stand out the power loss calculation using the nodal admittance matrix.  
\begin{example}
The following script creates a quadratic form using a random matrix $M\in\mathbb{R}^{2\times 2}$.  Notice how a function is created and used in Julia.
\begin{lstlisting}[frame=trBL]
function q(x)
    M = [3.5 4.2;8.1 5.0]
    return transpose(x)*M*x 
end
x0 = [1,1]
y0 = q(x0)
println(y0)
\end{lstlisting}
The operation \texttt{transpose(x)} can be replaced by $x'$ with the same result. However, it is important to be careful in the case of complex arrays where $x'$ transpose and conjugate. Therefore, it is recommended to use the command transpose to be certain we are only transposing without conjugating.
\end{example}

\begin{example}
The shape of a quadratic form may vary significantly according to the properties of its generating matrix.  Figure \ref{fig:plot_quadratic} shows two different quadratic forms: the first one generated by the matrix \eqref{eq:MMMM} and the second one by the following matrix: 
\begin{equation}
M = \begin{bmatrix} 0.6 & 0.2 \\ 0.2 & 0.8 \end{bmatrix}
\end{equation}
The second plot is a paraboloid, which is a convex function generated by a positive definite matrix (a concept that will be defined below). These plots can be generated in Julia as follows:
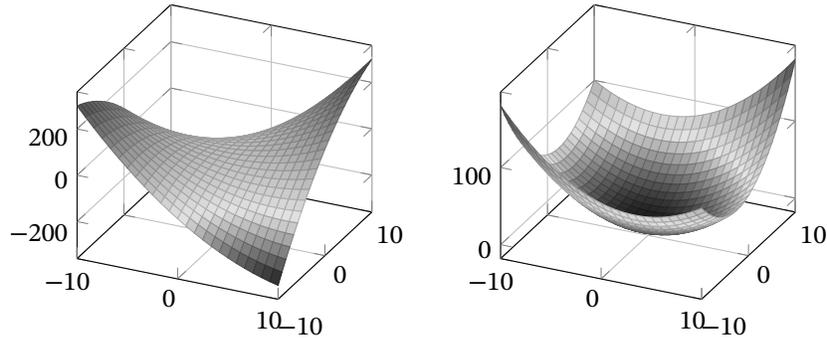
\begin{figure}[hbt]
    \centering
    \begin{tabular}{cc}
    \begin{tikzpicture} 
    \begin{axis}[grid=major, width=5.5cm, height=5.5cm]
	\addplot3[surf,domain=-10:10] {x^2+3*x*y-y^2};
	\end{axis}    
    \end{tikzpicture}  &
    \begin{tikzpicture}
    \begin{axis}[grid=major, width=5.5cm,height=5.5cm]
	\addplot3[surf,domain=-10:10] {0.6*x^2+0.4*x*y+0.8*y^2};
	\end{axis} 
    \end{tikzpicture}
    \end{tabular}
    \caption{Example of two different quadratic forms: left) $q(x_1,x_2)=x_1^2+3x_1x_2-x_2^2$, right) $q(x_1,x_2)=0.6x_1^2+0.4x_1x_2+0.8x_2^2$. The latter is a paraboloid defined by a quadratic form with a positive definite matrix. }
    \label{fig:plot_quadratic}
\end{figure}

\begin{lstlisting}[frame=trBL]
using Plots
xs = range(-10, stop=10, length=100)
ys = range(-10, stop=10, length=100)
M1 = [0.6 0.2;0.2 0.8]
f1(x1,x2) = [x1;x2]'*M1*[x1;x2] 
M2 = [1 3;0 -1]
f2(x1,x2) = [x1;x2]'*M2*[x1;x2] 
plt1 = surface(xs, ys, f1)
plt2 = contourf(xs, ys, f1)
plt3 = surface(xs, ys, f2)
plt4 = contourf(xs, ys, f2)
plot(plt1,plt2,plt3,plt4)
\end{lstlisting}
\end{example}

We say that $M$ is positive definite, represented by $M\succ 0$, if $q(x)>0$ for any $x\in\mathbb{R}^n$ different from zero.  Positive definite matrices are non-singular, meaning they have an inverse and its determinant is different from zero. On the other hand, a matrix is positive semidefinite, i.e., $M\succeq 0$, if $q(x)\geq 0,\; \forall x\neq 0$.  Positive definite and semidefinite matrices play a crucial role in many applications. For example, a quadratic form with a positive semidefinite matrix is convex.  This property is important in optimization problems since convexity allows guaranteeing global optimum and convergence of the algorithms \cite{9640361}. It also helps to demonstrate convergence in the fixed-point and newton methods for the power flow analysis \cite{yonewton}. 

\begin{example}
All eigenvalues of a positive definite matrix are positive. Hence, we can calculate the eigenvalues to determine if a matrix is positive definite. The following script calculates the eigenvalues \texttt{L} and eigenvectors \texttt{V} of a $2\times 2$ matrix $A$:
\begin{lstlisting}[frame=trBL]
using LinearAlgebra
A = [3 -2; -2 8]
(L,V) = eigen(A)
println("Eigenvalues \t=",L)
println("Eigenvectors\t=",V)
if minimum(L)>0 
    println("The matrix is positive definite")
end
\end{lstlisting}
This example was also an excuse to show the use of the \texttt{if} statement, which is intuitively defined in Julia.
\end{example}
We have the option to determine whether a matrix is positive definite by calculating its eigenvalues. However, certain matrices possess specific structural properties that guarantee positivity (or negativity) definiteness without the need to compute their eigenvalues. For instance, a matrix that is strictly diagonally dominant qualifies as positive definite. This type of matrix has the following property:
\begin{equation}
a_{kk} > \sum_{m\neq k} |a_{km}| \label{eq:diagonalmente_dominante}
\end{equation}
where $a_{kk}$ represents the diagonal entries of a matrix $A$ and, $a_{km}$ represents the non-diagonal entries $(m\neq k)$.  If the inequality is not strict, that is, if we replace $>$ by $\geq$ in \eqref{eq:diagonalmente_dominante}, then we say the matrix is weakly diagonal dominant.  The $Y_\text{bus}$ is the archetypical example in power systems.  Weakly diagonal matrices are positive definite, but this holds true only when they are also weakly chained (see \cite{weaklychained}). If they fail to meet this additional condition, then we can only guarantee that they are positive semidefinite. \index{diagonally dominant matrix} \index{weakly diagonal dominant matrix}

A positive definite matrix is non-singular and hence, the inverse is well-defined. Moreover, it can be represented by $A=C^\top C$ where $C$ is a triangular matrix named Cholesky factorization; $C$ is represented as $A^{1/2}$ since $A=C^\top C$ (not to confuse with a elementwise square root of $A$). This factorization can be used for fast solving of linear systems of the form $Ax=b$.

Like a quadratic form, a norm is another type of function that takes objects from $\mathbb{R}^n$ and returns a real number.  Intuitively, a norm is a way to measure the size of a vector.  For instance, the Euclidean norm, which is the most conventional of these functions defined as follows:
\begin{equation}
\left\| x \right\| = \sqrt{x^\top x}
\end{equation}
However, we can define a Q-norm from any quadratic function defined by a matrix $Q$ as long as this matrix is positive definite:
\begin{equation}
\left\| x \right\|_Q = \sqrt{x^\top Q x},\; Q\succ 0
\end{equation}
Similar to the Euclidean norm, $\left\| x\right\|_Q$ is positive for any $x$ except for the vector $x=0$. In addition, the triangular inequality holds, that is:
\begin{equation}
\left\| x+y\right\|\leq \left\| x\right\| + \left\| y\right\|
\end{equation}
Clearly, the numerical result of $\left\|x\right\|$ is different from $\left\|x\right\|_Q$; however, both functions are suitable ways to measure the size of a vector $x\in\mathbb{R}^n$. We may define a norm according to the application. Thus, $Q$ can be regarded as a weight factor.

Similarly, we can define a matrix norm as a manner to measure the size of a matrix. We induce a matrix norm given a vector norm as given below:
\begin{equation}
\left\| A\right\| = \underset{x\neq 0}{\sup} \frac{\left\| Ax\right\|}{\left\|x\right\|}
\label{eq:induced_norm}
\end{equation}
For any matrix $A$ and a vector $x$, we have the following property:
\begin{equation}
\left\| Ax\right\| \leq \left\| A\right\|\left\| x\right\|
\end{equation}

Finally, the Frobenius norm is another way to define a size of a matrix.  This norm is defined as follows:
\begin{equation}
\left\| A\right\|_F = \left(\sum_{i=1}^{n}\sum_{j=1}^{n} a_{ij}^2\right)^{1/2} = \sqrt{\operatorname{tr}(AA^\top)}
\label{eq:frobenius}
\end{equation}
Both, the induced norm given by \eqref{eq:induced_norm} and the Frobenius norm given by \eqref{eq:frobenius} satisfy the following property:
\begin{equation}
\left\| AB\right\| \leq \left\| A\right\|\left\| B\right\|
\end{equation}
In this case, we say that the matrix norm is sub-multiplicative. 
\begin{example}
Let us see the sub-multiplicative norm in practice.  The following script generates two random matrices $A$ and $B$; then, the norm of each one is calculated and compared:
\begin{lstlisting}[frame=trBL]
using LinearAlgebra
A = rand(2,2)
B = rand(2,2)
n_a = norm(A)
n_b = norm(B)
n_ab = norm(A*B)
println("norm(AB)=",n_ab, 
        " < norm(A)norm(B)=",n_a*n_b)
\end{lstlisting}
Frobenius norm is the default matrix-norm in Julia. Induced 2-norm can be obtained by the command \texttt{opnorm(A,2)} or simply \texttt{opnorm(A)}.
\end{example}

\section{Linear systems and sparse matrices}

Most of the power flows methods in for active distribution networks require solving a system of linear equations as follows:
\begin{equation}
Ax = b \label{eq:axb}
\end{equation}
where $A$ is a square matrix and $b$ is a column vector. Like most of the modern programming languages, Julia has methods for efficiently solving this type of problem. The most intuitive although less efficient way to solve \eqref{eq:axb} is, perhaps, by calculating the inverse of $A$ obtaining the value of $x$ as $x=A^{-1}b$.  However, this calculation can be performed more efficiently by using the backslash operator or a specific matrix factorization. The example below shows the concept. \index{linear systems}\index{$\setminus$ operator}\index{sparse matrix}\index{pseudoinverse}

\begin{example}\label{ex:Ax=b}
The following script creates a random system of linear equations that is solved in three different ways:
\begin{lstlisting}[frame=trBL]
using LinearAlgebra

n = 1000
A = rand(n,n)
b = rand(n,1)

@time xa = inv(A)*b

@time xb = A\b

L,U,p = lu(A)
@time begin
    y = L\b[p]
    xc = U\y
end
\end{lstlisting}
Here, \texttt{@time} is a macro for measuring calculation time.  First, a random matrix $A$ and a random vector $b$ are generated.  Next, the system is solved by three methods: direct calculation of the inverse, backslash operator and the $LU$ factorization.  The first two methods are intuitive. The third method is based on the following factorization:
\begin{equation}
PA = LU
\end{equation}
where $L$ is a lower triangular matrix, and $U$ an upper triangular matrix; $P$ is a permutation matrix.  Then, the linear system can be solved in two steps: first, the system $Ly=Pb$ is solved, by creating a new vector $y$ for the intermediate calculation; next, the system $Ux=y$ is solved to obtain the vector $x$. Notice that, Julia returns a permutation vector instead of a permutation matrix, so  the product $Pb$ can be efficiently calculated as \texttt{b[p]}. The student is invited to execute this code and compare perfomance. 
\end{example}

If the matrix is singular, we can still obtain a solution for the linear system \eqref{eq:axb} using the Moore-Penrose pseudoinverse. The pseudoinverse is represented as $A^+$ and calculated in Julia by the command \texttt{pinv()}.  Singularity of $A$ implies that the solution either does not exist or it is not unique. Thus, the pseudoinverse provides the least squares solution to the linear system, that is, a solution with a minimum Euclidean norm.  In this case, the value of $x$ is the best possible vector associated to the following optimization problem:
\begin{equation}
  \min \;  \left\|Ax-b\right\|^2 + \mu \left\|x\right\|^2
\end{equation}
where $\mu\rightarrow 0$.  Mathematical properties of the pseudoinverse can be studied in \cite{pseudoinverse}

\begin{example}
The script presented below is a simple calculation for a linear system with a singular matrix. 
\begin{lstlisting}[frame=trBL]
using LinearAlgebra
A = [1 3; 2 6]
b = [8;16]
x = pinv(A)*b
println(x)
println(A*x-b)
\end{lstlisting}
We can prove the matrix is singular by calculating its determinant.
\end{example}

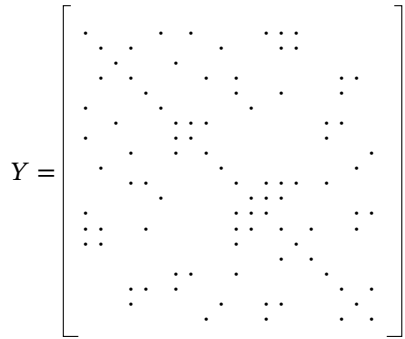
\begin{figure}[hbt]
\centering
\begin{tikzpicture}[x=1mm,y=-1mm]
    \node at (-5,20) {$Y=$};
    \draw (0,-2) -| +(-1,44) -- +(1,44);
    \draw (43,-2) -| +(1,44) -- +(-1,44);
\node at (2,2) {$ \cdot $};\node at (2,12) {$ \cdot $};\node at (2,16) {$ \cdot $};\node at (2,26) {$ \cdot $};\node at (2,28) {$ \cdot $};\node at (2,30) {$ \cdot $};\node at (4,4) {$ \cdot $};\node at (4,8) {$ \cdot $};\node at (4,20) {$ \cdot $};\node at (4,28) {$ \cdot $};\node at (4,30) {$ \cdot $};\node at (6,6) {$ \cdot $};\node at (6,14) {$ \cdot $};\node at (8,4) {$ \cdot $};\node at (8,8) {$ \cdot $};\node at (8,18) {$ \cdot $};\node at (8,22) {$ \cdot $};\node at (8,36) {$ \cdot $};\node at (8,38) {$ \cdot $};\node at (10,10) {$ \cdot $};\node at (10,22) {$ \cdot $};\node at (10,28) {$ \cdot $};\node at (10,36) {$ \cdot $};\node at (12,2) {$ \cdot $};\node at (12,12) {$ \cdot $};\node at (12,24) {$ \cdot $};\node at (14,6) {$ \cdot $};\node at (14,14) {$ \cdot $};\node at (14,16) {$ \cdot $};\node at (14,18) {$ \cdot $};\node at (14,34) {$ \cdot $};\node at (14,36) {$ \cdot $};\node at (16,2) {$ \cdot $};\node at (16,14) {$ \cdot $};\node at (16,16) {$ \cdot $};\node at (16,34) {$ \cdot $};\node at (18,8) {$ \cdot $};\node at (18,14) {$ \cdot $};\node at (18,18) {$ \cdot $};\node at (18,40) {$ \cdot $};\node at (20,4) {$ \cdot $};\node at (20,20) {$ \cdot $};\node at (20,38) {$ \cdot $};\node at (22,8) {$ \cdot $};\node at (22,10) {$ \cdot $};\node at (22,22) {$ \cdot $};\node at (22,26) {$ \cdot $};\node at (22,28) {$ \cdot $};\node at (22,30) {$ \cdot $};\node at (22,34) {$ \cdot $};\node at (24,12) {$ \cdot $};\node at (24,24) {$ \cdot $};\node at (24,26) {$ \cdot $};\node at (24,28) {$ \cdot $};\node at (26,2) {$ \cdot $};\node at (26,22) {$ \cdot $};\node at (26,24) {$ \cdot $};\node at (26,26) {$ \cdot $};\node at (26,38) {$ \cdot $};\node at (26,40) {$ \cdot $};\node at (28,2) {$ \cdot $};\node at (28,4) {$ \cdot $};\node at (28,10) {$ \cdot $};\node at (28,22) {$ \cdot $};\node at (28,24) {$ \cdot $};\node at (28,28) {$ \cdot $};\node at (28,32) {$ \cdot $};\node at (28,38) {$ \cdot $};\node at (30,2) {$ \cdot $};\node at (30,4) {$ \cdot $};\node at (30,22) {$ \cdot $};\node at (30,30) {$ \cdot $};\node at (32,28) {$ \cdot $};\node at (32,32) {$ \cdot $};\node at (34,14) {$ \cdot $};\node at (34,16) {$ \cdot $};\node at (34,22) {$ \cdot $};\node at (34,34) {$ \cdot $};\node at (36,8) {$ \cdot $};\node at (36,10) {$ \cdot $};\node at (36,14) {$ \cdot $};\node at (36,36) {$ \cdot $};\node at (36,40) {$ \cdot $};\node at (38,8) {$ \cdot $};\node at (38,20) {$ \cdot $};\node at (38,26) {$ \cdot $};\node at (38,28) {$ \cdot $};\node at (38,38) {$ \cdot $};\node at (40,18) {$ \cdot $};\node at (40,26) {$ \cdot $};\node at (40,36) {$ \cdot $};\node at (40,40) {$ \cdot $};
    
\end{tikzpicture}
\caption{Visualization of a sparse matrix.  Dots represent entries different from zero.}
\label{fig:spartsematrix}
\end{figure}

On the other hand, most of the admittance matrices in power systems are highly sparse.  This type of matrix has a characteristic structure where most of its entries are zero, as depicted in Figure \ref{fig:spartsematrix}.  In this case, we can efficiently store and handle the matrix using the package \texttt{SparseArrays}. Calculations with a sparse matrix may be done efficiently since it has fewer actual data points compared to dense matrices. When performing calculations on sparse matrices, algorithms can take advantage of this sparsity by skipping operations involving zero elements. This reduces the overall number of computations needed, leading to faster computation times. Using sparse arrays is simple, as demonstrated in the example below.

\begin{example}
The following script creates a sparse matrix $A=M+D$, where each entry in $M$ is generated by rounding a random number between 0 and 0.6, and, $D$ is a random diagonal matrix resulting in a sparse matrix. The rest of the code is intuitive.
\begin{lstlisting}[frame=trBL]
using LinearAlgebra
using SparseArrays
n = 10
M = round.(0.6*rand(n,n))
D = Diagonal(rand(n))
A = sparse(M+D)
display(A)
\end{lstlisting}

\end{example}

\section{Data frames}

Power distribution networks may have thousands of nodes and lines, which imply large amount of data. Therefore, it is required to efficiently store and handle information through data frames. This data structure organizes the information into a 2-dimensional table of rows and columns, similarly to a spreadsheet. The package \texttt{DataFrames} allows intuitively handling this type of structure, as exemplified below. 

\begin{example}
The information of four substations is stored in three vectors that contain the identification name, the voltage level in kV and, whether it is connected to the system. This information is registered in vectors \texttt{A,B,C}.  Then, a data frame is created.  Finally, the information of the data frame is printed, first as the entire data frame and next individually. 
\begin{lstlisting}[frame=trBL]
using DataFrames
A = ["Electron platyrhynchum",
     "Guacamaya",
     "Motmot",
     "Trochilidae"]
B = [500;500;115;230]
C = [true,true,true,false]
df = DataFrame(Name=A,
               Voltage_level=B, 
               Connected=C)
println(df)
# accesing individual data
println(df.Name[1])
println(df.Voltage_level[3])
using Plots
bar(df.Name,df.Voltage_level)
\end{lstlisting}
Notice the data-frame can store different type of data (e.g., text, numbers and boolean) in the same structure.  The last line combines the package \texttt{DataFrame} with the package \texttt{Plots} to obtain a bar plot. The student is invited to execute the code above to see the way as a data frame organizes the information.
\end{example}

It is common in power systems applications to storage the information of the test systems using a comma separated values or a CSV file.  This is a plain text file that stores data by delimiting data entries with commas. Julia is able to read said file, transforming into a data frame object. For that, we require the packages \texttt{CSV} and, \texttt{DataFrames}.  Let us consider the following practical example: 

\begin{example}\label{ex:tomaruncvs}
Let us consider a CSV file with the data presented in Table \ref{tab:data_cigre}, which corresponds to the CIGRÉ benchmark for medium voltage distribution network \cite{cigre_microgrid}. In this table, the values of R and X are in Ohm/km, C in nF/km and Length in km.  We named this file as \texttt{CigreMicrogrid.csv}, which is loaded with the following script, saved in the same folder as the script:

\begin{table}[hbt]
\centering
\caption{Data for the CIGRÉ medium voltage distribution networks taken from \cite{cigre_microgrid}}
\label{tab:data_cigre}
\begin{tabular}{rrrrrr}
\toprule
From & To &R & X & C&    Length \\
\midrule
1& 2&   0.579& 0.367& 993& 2.82 \\
2& 3&   0.164& 0.113& 413& 4.42 \\
3& 4&   0.262& 0.121& 405& 0.61 \\
4& 5&   0.354& 0.129& 285& 0.56 \\
5& 6&   0.336& 0.126& 343& 1.54 \\
6& 7&   0.256& 0.130& 235& 0.24 \\
7& 8&   0.294& 0.123& 350& 1.67 \\
8& 9&   0.339& 0.130& 273& 0.32 \\
9&  10& 0.399& 0.133& 302& 0.77 \\
10& 11& 0.367& 0.133& 285& 0.33 \\ 
3&  8&  0.172& 0.115& 411& 1.30 \\
11& 4&  0.423& 0.134& 310& 0.49 \\
\bottomrule
\end{tabular}
\end{table}

\begin{lstlisting}[frame=trBL]
using CSV
using DataFrames
microgrid = DataFrame(CSV.File("CigreMicrogrid.csv"))
print(microgrid)
\end{lstlisting}

Entries in a data frame can be obtained by the index number; for example, \texttt{microgrid[1,3]} returns $0.579$ corresponding to the resistance of the first line section which connects nodes 1 and 3.  The same value can be obtained as \texttt{microgrid.R[1]}. In fact, all resistances can be obtained as an array given by \texttt{microgrid.R[:]}. However, printing a dataframe results in neater output than printing an array. We finish this example by visualizing the values of R as a bar plot, as given below:
\begin{lstlisting}[frame=trBL]
using Plots
bar(microgrid.R,label="Resistance",
xlabel="Node",ylabel="Ohm/km")
\end{lstlisting}
\end{example}

\section{Good practices}

Writing a general purpose software is different from writing a scientific software. However, there are some concepts from software engineering that may be useful to obtain a clean and readable script. Here, we present some recommendations inspired in
\cite{practices_software}:
\begin{description}
\item[Planning:] think about the architecture of the script before writing the code.  Decide how many functions are required and how the information will flow among these functions. Drawing on a piece of paper may help to see the entire architecture, thinking about the future and how the script may growth. 
\item[Writting:] place a brief explanatory comment at the start of every program and function. Give functions and variables meaningful names.  Julia allows including Greek letters and tildes as variables.  This is a very useful feature that should be used to the fullest. Include a README file with details about the project. Markdown is an excellent tool for this purpose.
\item[Organize the input data:] It is important to properly organize CSV or Excel databases before running any code.  Organize the data as simple as possible, similarly to a Dataframe. Choose meaningful names for each variable. Use a single word without spaces for each variable.  For example,  \texttt{Rpu} instead of \texttt{R} or \texttt{R pu} to represent resistance when it is given in per unit. See for example \cite{data_organization} for recommendations about data organization.
\item[Create a toy-model:] do not pretend to solve a large test system without solving first a small ``toy-model''. For instance, you may consider first a small DC distribution network with three or four nodes, or an AC system without reactive power, before considering the full problem.
\item[Decompose programs into functions:] Breaking code into functions allows creating modular components that may be used in further applications. For example, it is useful to create functions for calculating the $Y_\text{bus}$ which is used in almost every power system application.  These functions should be small and easy to read. The Julia compiler can often better optimize individual functions compared to a single large script.
\end{description}

\section{Further readings}

There is never enough linear algebra for a power systems engineer. Several books, web pages and even videos address to this subject. For example, the work of \cite{linear_algebra_julia} presents a complete course of linear algebra with applications in Julia.  A more formal treatment can be found in books like  \cite{linear_algebra_don_wright}, \cite{hubbard}, and \cite{horn_johnson_1985}. It is also important to learn how to transform mathematical concepts into useful algorithms (this is one of the primary motivations of this book). However, writing clean code is more an art than a science. There is plenty of material available in books and on the Internet. For instance, \cite{martin2009clean} presents a comprehensive introduction to software engineering. Nevertheless, it is essential to remember that power systems' computation is a field with particular characteristics. What works in an accounting software may not be relevant in power systems applications. The experience helps to develop better scripts, so, it is recommended to try different approaches until finding the best practices.

\bibliographystyle{apalike}
\bibliography{bibliografia}

\end{document}